\documentclass[11pt]{article}

\pdfoutput=1

\RequirePackage{my_packages}
\RequirePackage{my_theorems}
\RequirePackage{my_macros}

\usepackage[backend = biber,        
            language = english ,
            style = alphabetic ,   
            giveninits = true,
            isbn = false,
            url = false,
            doi = false,
            sorting = nyt,
            maxnames = 6,
            backref=true
            ]{biblatex}

\renewbibmacro{in:}{%
  \ifentrytype{article}{}{%
  \printtext{\bibstring{in}\intitlepunct}}}
\AtEveryBibitem{%
  \clearfield{note}%
  \clearfield{edition}%
  \clearlist{language}
}

\DeclareNameAlias{labelname}{given-family}

\title{Decidability of the Brinkmann Problems\\for endomorphisms of the free group}
\author{André Carvalho$^*$ \& Jordi Delgado$^\dag$}
\date{\vspace{-7pt}
    $^*$Center for Mathematics and Applications (NOVA Math), NOVA FCT\\%
    $^\dag$Departament de Matemàtiques, Universitat Politècnica de Catalunya\\[3ex]%
    \today
}

\newcommand{\Addresses}{{
  \bigskip
  \footnotesize

  André Carvalho\\\nopagebreak
  \textsc{Center for Mathematics and Applications (NOVA Math), NOVA FCT}\\\nopagebreak
  \url{andrecruzcarvalho@gmail.com}

  \medskip

  Jordi Delgado\\\nopagebreak
  \textsc{Departament de Matemàtiques, Universitat Politècnica de Catalunya}\\\nopagebreak
  \url{jorge.delgado@upc.edu }

}}
\hypersetup{final}  

\begin{document}

\maketitle


\begin{abstract}
\noindent 
Building on the work of Brinkmann and Logan, we show that both the Brinkmann Problem
and the Brinkmann Conjugacy Problem
are decidable for endomorphisms of the free group $\Fn$.
\end{abstract}

 \bigskip



\vspace{15pt}

Let $\NN$ denote the set of natural numbers (including $0$), and let $\Fn$ denote the free group of finite rank $n$.
In this note, we elaborate on previous work of Brinkmann 
and Logan
to prove the result below.

\begin{thm*}
    Given two elements $u,v \in \Fn$ and an endomorphism $\varphi \in \End(\Fn)$, it is algorithmically decidable whether:
    \begin{enumerate}[ind]
        \item there exists some $k \in \NN$ such that $(u) \varphi^k = v$;
        \item there exists some $k \in \NN$ such that $(u)\varphi^k$ is conjugate to $v$. \qed
    \end{enumerate}  
\end{thm*}

That is, both the \defin{Brinkmann problem} (\BrP) and the Brinkmann Conjugacy Problem (\BrCP) are decidable for endomorphisms of free groups. This result generalizes the analogous claims on automorphisms and monomorphims of the free group, proved by Brinkmann \cite{brinkmann_detecting_2010} and Logan \cite{logan_conjugacy_2023} respectively.

\section{Background}
The study of algorithmic problems for groups dates back to the early twentieth century, when Max Dehn proposed three seminal decision problems (namely, the \defin{word problem}, the \defin{conjugacy problem}, and the \defin{isomorphism problem}) in \cite{dehn_uber_1911}.
The scope of the first two of them is a group $G$ given by a finite presentation $\pres{X}{R}$:

The \defin{Word problem} for $G$, denoted by $\WP(G)$,  consists in deciding, given two words $u,v \in (X^{\pm})^*$, whether they represent the same element in $G$. 

The \defin{Conjugacy problem} for $G$, denoted by $\CP(G)$,  consists in deciding , given  $u,v \in (X^{\pm})^*$, whether they represent conjugate elements in $G$.

The problems of this kind (admitting as inputs strings $w \in \Sigma^*$ in some finite alphabet $\Sigma$, and asking whether $w$ belongs to a certain subset $S \subseteq \Sigma^*$) are called \defin{decision problems}. A decision problem is said to be \defin{decidable}
if there exists an algorithm (formally, a Turing machine) which, on every input $w$, outputs \yep\ if $w \in S$, and \nop\ if $w \notin S$.


It is not difficult to see using Tietze transformations that the outcome of all the problems appearing in this paper does not depend on the chosen presentation. Accordingly, it is common to use the elements of the group and the words representing them (in some presentation) interchangeably, provided that the meaning remains clear from the context.


Natural variations of the previous decision problems soon arose, for example the \defin{Subgroup Membership Problem}, denoted by $\MP(G)$, consists in deciding, given finitely many words $w_0,w_1,\ldots,w_k$ in $X^{\pm}$, whether the group element represented by $w_0$ belongs to the subgroup generated in~$G$ by~$w_1,\ldots,w_k$. 


Recall that, for an arbitrary group $G$ and $g_1,g_2\in G$, we
say that $g_1$ is \defin{conjugate} to $g_2$ (in $G$),
and we write $g_1 \conj g_2$  (or $g_1 \conj_G g_2$),
if there exists some element $h \in G$ such that~$g_1 = h^{-1} g_2 h$.
It is routine to check that conjugacy is an equivalence relation in $G$ whose quotient set, denoted by $G/{\conj}$, does not necessarily inherit the group operation.
Also, 
if $S \subseteq G$, we write $g \conj S$ to mean that $g \conj h$ for some $h \in S$.

In the 1930s, endomorphisms were added to the picture through  the now called \emph{Whitehead's orbit problems}: given $u,v \in (X^{\pm})^*$, decide whether there exists some automorphism $\varphi \in \Aut(G)$ such that $(u)\varphi =_G v$  (\resp whether $(u)\varphi \conj_G v$). If $G = \pres{X}{R}$, an endomorphism $\varphi \in \End(G)$ is given as an input through the set $(X)\varphi$ of images under $\varphi$ of the generators for $G$.

Although all of the mentioned problems
were famously shown to be algorithmically undecidable for a generic finitely presented group in the mid-20th century (see \cite{novikov_algorithmic_1955,boone_word_1958}), 
in the context of finitely generated free groups $\Fn = \pres{x_1,\ldots x_n}{-}$, they are well-known to be  decidable: $\WP(\Fn)$ and $\CP(\Fn)$ using elementary arguments involving reduced and cyclically reduced words; $\MP(\Fn)$ using, for example, the nice theory of Stallings automata (see~\cite{stallings_topology_1983,Delgado_Stallings_automata_GAGTA2024}); and Whitehead orbit problems using the now classical \defin{peak-reduction} technique (see~\cite{whitehead_equivalent_1936}).

\section{Map orbits and the Brinkmann Problems}

At the beginning of the 21st century, 
equipped with
the recently developed theory of relative train track maps  (\cite{bestvina_train_1992,bestvina_train-tracks_1995}),    \citeauthor{brinkmann_detecting_2010} \cite{brinkmann_detecting_2010}  
addressed two cyclic variations of the orbit problem for free groups.
Concretely, he proved that, given two elements $u,v \in \Fn$, and an automorphism $\varphi \in \Aut(\Fn)$, it is always possible to algorithmically decide 
whether there exists some
exponent $k$
such that~$u\varphi^k = v$ (\resp whether $u\varphi^k$ is conjugate to $v$).

Note the ambiguity in the term 'exponent' in the claim above. Is it meant to be an integer? a natural number? does it include zero? Below we specify the different variants that arise, we clarify the relation between them, and
we prove that all of them are algorithmically equivalent 
in Brinkmann's original context.

\begin{defn}
    Let $G$ be a set, let $u \in G$, and let $\varphi\colon G \to G$ be a map. Then,
    \begin{itemize}
    \item the \defin{strict orbit} of $u$ under $\varphi$ 
    is $\orb_{\varphi}^{+}(u) = 
        \set{u\varphi^k \st k > 0}$;
    \item the \defin{(standard) orbit} of $u$ under $\varphi$ 
    is $\orb_{\varphi}(u) = 
        \set{u\varphi^k \st k \in \NN}
        $;
    \item the \defin{strict symmetric orbit}  of $u$ under $\varphi$ is $\corb_{\varphi}^{\,+}(u) = 
         \set{v \in G \st v\in \orb^{+}_{\varphi}(u) \text{ or } u\in \orb^{+}_{\varphi}(v) }$;
    \item the \defin{symmetric orbit} of $u$ under $\varphi$ is $\corb_{\varphi}(u) =
        \set{v \in G \st v\in \orb_{\varphi}(u) \text{ or } u\in \orb_{\varphi}(v) }$.    
    \end{itemize}
    \end{defn}

Orbits under $\varphi$ are usually referred as $\varphi$-orbits. Note that
the symmetric $\varphi$-orbit of $u$ is the union of all the $\varphi$-orbits containing $u$. Also recall that if $\varphi
    $ is bijective then
    $\corb_{\varphi}^{\,+} (u) = \orb_{\varphi}^{+} (u) \cup \orb_{\varphi^{\text{-}1}}^{+} (u)
    = \set{u \varphi^k \st k \in \ZZ \setmin \set{0}}$; 
    and
    $\corb_{\varphi} (u) = \orb_{\varphi} (u) \cup \orb_{\varphi^{\text{-}1}} (u)$.
 
\begin{rem}
    If $\varphi \colon G \to G$ is a homomorphism of groups, then $\widetilde{\varphi} \colon G/{\conj} \to G/{\conj}$, $[u] \mapsto [u \varphi]$ is a well defined map, 
and hence it makes sense to consider any of the previous versions of orbits under $\widetilde{\varphi}$ (we call them \defin{orbits of $\varphi$ modulo conjugacy}).
\end{rem}

\begin{defn}
An element $u \in G$ is said to be \defin{periodic} \wrt $\varphi$ (or \defin{$\varphi$-periodic}) if it belongs to its own strict orbit; that is, if $u = u \varphi^k$ for some $k > 0$. If $u$ is $\varphi$-periodic, the minimum  integer $p>0$ such that $u \varphi^p = u$ is called the \defin{$\varphi$-period} of $u$. The set of $\varphi$-periodic elements in $G$ is denoted by~$\Per(\varphi)$.
Accordingly, if $G$ is a group, we say that $u \in G$ is \defin{conjugate-periodic} \wrt $\varphi$ 
(or \defin{$\varphi$-conjugate-periodic}) 
if $[u]\in G/{\conj}$ is $\widetilde{\varphi}$-periodic
(that is, if there exists some $k>0$ such that $u \varphi^k \conj u$). The \defin{conjugate $\varphi$-period} of $u$ is the $\widetilde{\varphi}$-period of $[u]$, and the set of conjugate-periodic elements under $\varphi$ is denoted by $\CPer(\varphi)$.
\end{defn}

The algorithmic problem(s) consisting in deciding whether a given element is periodic \wrt certain given maps are called \defin{Orbit Periodicity Problems}. 

\begin{defn}
    Let $G$ be a group and let $\set{\id_G} \subseteq \mathcal{T} \subseteq \End(G)$. Then,
    \begin{itemize}
        \item the \defin{Orbit Periodicity Problem} 
for $G$ \wrt $\T$, denoted by $\OPP_{\T}(G)$, consists in deciding, 
given $u \in G$ and $\varphi \in \T$, whether $u$ is $\varphi$-periodic.
        \item the \defin{Orbit Conjugate-Periodicity Problem} 
for $G$ \wrt $\T$, denoted by $\OCPP_{\T}(G)$, consists in deciding, 
given $u \in G$ and $\varphi \in \T$, whether $u$ is $\varphi$-conjugate-periodic.
    \end{itemize}
\end{defn}

\begin{rem}
If $\OPP_{\T}(G)$ and $\WP(G)$ (\resp $\OCPP_{\T}(G)$ and $\CP(G)$) are decidable, then the orbit and period of any $\varphi$-periodic (\resp $\varphi$-conjugate-periodic) element is computable just by inspection.  
\end{rem}


We agglutinate under the name of \defin{Brinkmann Problems} the problems consisting in deciding whether two elements belong to the same `orbit' under a map.
\begin{defn} \label{defn: Brinkmann}
    Let $G$ be a group and let $\set{\id_G} \subseteq \mathcal{T} \subseteq \End(G)$. Then,
    \begin{itemize}
        \item the \defin{(standard) Brinkmann Problem} 
for $G$ \wrt $\T$, denoted by $\BrP_{\T}(G)$, consists in deciding, given $u,v \in G$ and $\varphi \in \T$, whether $v \in \orb_{\varphi}(u)$.
        \item the \defin{(standard) Brinkmann Conjugacy Problem} 
for $G$ \wrt $\T$, denoted by $\BrCP_{\T}(G)$, consists in deciding, given $u,v \in G$ and $\varphi \in \T$, whether $[v] \in \orb_{\widetilde{\varphi}}([u])$.
    \end{itemize}
\end{defn}
The natural variants for
 strict, symmetric, or strict-symmetric 
orbits are denoted by
$\BrP_{\T}^{+}(G)$, 
$\overline{\BrP}_{\T}(G)$ and 
$\overline{\BrP}_{\T}^{\,+}(G)$
for the corresponding Brinkmann Problems; and
$\BrCP_{\T}^{+}(G)$, 
$\overline{\BrCP}_{\T}(G)$ and 
$\overline{\BrCP}_{\T}^{\,+}(G)$
 for the corresponding Brinkmann Conjugacy Problems.

We abbreviate $\BrPa(G) = \BrP_{\!\Aut(G)}(G)$, $\BrPm(G) = \BrP_{\Mon(G)}(G)$, and $\BrPe(G) = \BrP_{\End(G)}(G)$, and similarly for the other variants.

\begin{rem} \label{rem: BP(id) = WP}
Note that $\BrP_{\mathrm{id}_G}(G)  = \WP(G)$ and   $\BrCP_{\mathrm{id}_G}(G) =\CP(G)$. 
Hence, it is clear that
$\WP(G) \preceq \BrPa(G) \preceq \BrPm(G) \preceq \BrPe(G)$ and 
$\CP(G) \preceq \BrCPa(G) \preceq \BrCPm(G) \preceq \BrCPe(G)$.\footnote{\label{footnote}We write  $\mathsf{P} \preceq \mathsf{Q}$ to express that
 $P$ can be solved using a potential algorithm solving $Q$.}

(Following the same reasoning, analogous inclusions can be obtained for each of the variants of the Brinkmann Problem in the paragraph following \Cref{defn: Brinkmann}.)
\end{rem}

\begin{rem} \label{rem: yes outputs}
The decidability of (any of) the Brinkmann's problems
allows to compute a witness in case it exists: if the answer of $\BrP_{\T}(G)$ on input $(u,v,\varphi)$ is \yep, then it is enough to keep enumerating the successive images $((u)\varphi^n)_{n \geq 0}$ and, in parallel, enumerate the consequences of the relators until reaching a guaranteed word of the form $v^{-1} (u)\varphi^k$, providing a witness~$k \in \NN$.
In fact,
since $\WP(G) \preceq \BrPa(G)$ (\resp $\CP(G) \preceq \BrCPa(G)$), we can easily compute the full set of witnesses for Brinkmann's problems (see \cite{carvalho_free-abelian_2024} for details). 
\end{rem}

The relation between the different variants of the Brinkmann problem is clarified below.

\begin{prop} \label{prop: equivalence}
     Let $G$ be a group, and let  $\set{\id_G} \subseteq \mathcal{T} \subseteq \End(G)$. Then, the following statements are equivalent:
    \begin{enumerate}[dep]
     \item $\BrP_{\T}(G)$ is decidable;
    \item $\BrP_{\T}^{+}(G)$ is decidable;
    \item $\overline{\BrP}_{\T}^{\,+}(G)$ is decidable;
    \item both $\overline{\BrP}_{\T}(G)$ and $\OPP_{\T}(G)$ are decidable.
    \end{enumerate}
\end{prop}

     


\begin{proof}
We recall that in any of the four cases, the outcomes of the inputs of the form $(u,v,\id_G)$ correspond to $\WP(G)$, which therefore can be assumed to be decidable throughout the proof.

To see that (a)$\Imp$(b), note that the output of
$\BrP_{\T}^{+}(G)$ on input $(u,v,\varphi)$ is the same as the output of $\BrP_{\T}(G)$ on input $(u\varphi,v,\varphi)$. 

The implication (b)$\Imp$(c) is also obvious since deciding whether $v \in \overline{\orb}_{\varphi}^{\,+}(u)$ immediately reduces to two instances of $\BrP^{+}(G)$.

To see that (c)$\Imp$(d), first note that, on inputs of the form $(u,u,\varphi)$, problem $\overline{\BrP}_{\T}^{\,+}(G)$ corresponds exactly to $\OPP_{\T}(G)$. On the other hand, the decidability of $\cBrP_{\T}(G)$ follows clearly from those of $\WP(G)$ (which is implied by $\cBrP_{\T}^{\,+}(G)$) and $\cBrP_{\T}^{\,+}(G)$: given an input $(u,v,\varphi)$, use $\WP(G)$ to check whether $u =_G v$, and $\cBrP_{\T}^{\,+}(G)$ to check whether $v \in \corb^{\,+}(u)$. The output of $\cBrP_{\T}(G)$ is \nop\ if both answers are negative, and \yep\ otherwise.

Finally, to see that (d)$\Imp$(a), let us assume that both $\cBrP_{\T}(G)$ and $\OPP_{\T}(G)$ are decidable, and consider an input $(u,v,\varphi)$ for $\BrP_{\T}(G)$.
We start giving $(u,v,\varphi)$ as an input for $\cBrP_{\T}(G)$. If the answer is \nop\ then neither $u \varphi^k = v$ nor $u \varphi^k = v$ (for $k \geq 0$) is possible, and the answer to $\BrP_{\T}(G)$ on the same input is \nop\ as well.

So, let us finally assume that $\cBrP_{\T}(G)$ answers \yep\ on $(u,v,\varphi)$. This means that:
\begin{equation} \label{eq: 2 cond}
   \text{either $v \in \orb_{\varphi}(u)$ \,or\, $u \in \orb_{\varphi}(v)$ \ (or both).} 
\end{equation}

Since at least one of the conditions in \eqref{eq: 2 cond} holds, we can start enumerating in parallel the elements in $\orb_{\varphi}(u)$ and $\orb_{\varphi}(v)$ until either $v$ appears in $\orb_{\varphi}(u)$ or $u$ appears in $\orb_{\varphi}(v)$. In the first case $v \in \orb_{\varphi}(u)$ and hence the answer to $\BrP_{\T}(G)$ on input $(u,v,\varphi)$ is \yep. Otherwise (if $u \in \orb_{\varphi}(v)$), note that $v \in \orb_{\varphi}(u)$ (and hence the answer to $\BrP_{\T}(G)$ is \yep) if and only if $v$ is periodic, which we can check using $\OPP_{\T}(G)$. So, we have decided $\BrP_{\T}(G)$ for every possible input $(u,v,\varphi)$ and the proof of (d)$\Imp$(a) is complete.
\end{proof}

An analogous result can be obtained for for the conjugacy version of the problems (replacing equality by conjugacy and $\WP$ by $\CP$ in the proof).

\begin{prop} \label{prop: equivalence conj}
     Let $G$ be a group, and let  $\set{\id_G} \subseteq \mathcal{T} \subseteq \End(G)$. Then, the following statements are equivalent:
    \begin{enumerate}[dep]
     \item  $\BrCP_{\T}(G)$ is decidable;
    \item $\BrCP_{\T}^{+}(G)$ is decidable;
    \item $\overline{\BrCP}_{\T}^{\,+}(G)$ is decidable;
    \item both $\overline{\BrCP}_{\T}(G)$ and $\OCPP_{\T}(G)$ are decidable.
    \end{enumerate}
\end{prop}

That is, from a decidability perspective, we only have two different variants of Brinkmann orbit problems; namely $\BrP_{\T}(G)$ and $\overline{\BrP}_{\T}(G)$ (\resp $\BrCP_{\T}(G)$ and $\overline{\BrCP}_{\T}(G)$).

\section{Brinkmann Problems for the free group}
According to  \Cref{prop: equivalence,prop: equivalence conj}, all the equality (\resp conjugacy) variants of the Brinkmann Problem are equivalent in any ambient group with solvable $\OPP$ (\resp $\OCPP$). It follows from
\Cref{thm: Brinkmann Fn} that this is precisely the case when the ambient is a finitely generated free group (and, in particular, in the context of $\Aut(\Fn)$ originally considered by Brinkmann in \cite{brinkmann_detecting_2010}). Hence, we will no longer distinguish between the standard and symmetric forms of the Brinkmann Problem throughout this section.




\R{
}

The initial work of Brinkmann (proving the decidability of $\BrPa(\Fn)$ and $\BrCPa(\Fn)$) found natural continuation in \cite{logan_conjugacy_2023}, where Logan uses Mutanguha's recent techniques extending to endomorphisms the computability of the fixed subgroups  of~$\Fn$ (see \cite{mutanguha_constructing_2022}) to prove the decidability of \CP\ for ascending HNN-extensions of free groups. On his way to this result, Logan considers and solves several variants of the Brinkmann Problem
including $\BrPm(\Fn)$ and $\BrCPm(\Fn)$, and two-sided versions of both problems. 
While not answering the general question for endomorphisms, in \cite[Section 5]{logan_conjugacy_2023}, Logan proves two results for endomorphisms that we will use to extend his analysis to reach the decidability of the general endomorphism cases.

\begin{thm} \label{thm: Brinkmann Fn}
    The Brinkmann Problems $\BrPe(\Fn)$  and $\BrCPe(\Fn)$ are algorithmically decidable.
\end{thm}

\begin{proof}
    Let $u,v \in \Fn$ and $\varphi \in \End(\Fn)$ be our input (given in terms of a basis $\{x_1,\ldots,x_n\}$ of~$\Fn$; in particular $\varphi$ is given by the $n$-tuple of words $v_1 = (x_1)\varphi,\ldots,v_n = (x_n) \varphi$).
    
In order to 
decide whether there exists some $k \in \NN$ such that $u \varphi^k = v$ we start using the decidability of $\WP(\Fn)$ to check whether $u = v$; if so, the answer to $\BrPe(\Fn)$ is \yep\ (with $k = 0$). 
Otherwise,
we can use the decidability of $\MP(\Fn)$ to check whether $v \in \im(\varphi) = \gen{v_1,\ldots,v_n}$. 
If not, 
then the answer to $\BrPe(\Fn)$ is \nop\ as well. Finally,
if $v \in \im(\varphi)$, then we can easily compute (either by brute force, or more efficiently using Stallings automata, see \eg \cite{Delgado_Stallings_automata_GAGTA2024}) an expression $v = w(v_1,\ldots,v_n)$  (of $v$ in terms of $v_1,\ldots,v_n$) and, from it, a preimage $v'= w(x_1,\ldots,x_n) \in \varphi^{-1}(v)$. Now it is clear that, for~$k\geq 1$, $u \varphi^k = v$ if and only if $u \varphi^{k-1} \in v'\ker(\varphi)$, and we can use the algorithm from Theorem 5.5 in \cite{logan_conjugacy_2023} (with inputs $u,v'$, and $\varphi$) to decide about the existence of such an integer~$k$. Hence, $\BrPe(\Fn)$ is decidable, as claimed.


For the second claim, first note that Theorem 5.2 in \cite{logan_conjugacy_2023} provides an algorithm~$\mathfrak{A}$ that, on input $u,v \in \Fn$ and $\varphi \in \End(\Fn)$, decides whether there exists some $k \in \NN$ such that $u \varphi^k \conj v \ker(\varphi)$, and outputs one such $k$ if it exists. 

In order prove the decidability of $\BrCPe(\Fn)$
we use the previous algorithm~$\mathfrak{A}$ in a combined manner.
On one side we run $\mathfrak{A}$ on input $(u,v,\varphi)$. If $\mathfrak{A}(u,v,\varphi)$ answers \nop, it means that for every $k \in \NN$, 
$u\varphi^k \nsim v \in v \ker(\varphi)$,
and the answer to $\BrCPe(u,v,\varphi)$ is also \nop. Otherwise, $\mathfrak{A}(u,v,\varphi)$ outputs an $k \in \NN$ such that $u \varphi^k \conj v \ker(\varphi)$. It is important to note that this means that $u \varphi^{k+1} \conj v\varphi$.

On the other hand, we run $\mathfrak{A}$ on input $(v\varphi,v,\varphi)$. Now, if $\mathfrak{A}(v\varphi,v,\varphi)$ answers \nop\ then $v \varphi^{p} \nsim v \ker(\varphi)$ for every $p\geq 1$. But since $u \varphi^{k+1} \conj v\varphi$, this also means that for every $p \geq 1$,
$
u \varphi^{k+p}
\conj
v \varphi^p
\nconj
v \ker(\varphi)
$,
and, in particular, that $u \varphi^{k+p} \nconj v$, for every $p \geq 1$. That is, in this case there are only finitely many candidates and it is enough to use $\CP(\Fn)$ to check whether $u \varphi^j \conj v$ for some $j \in\{0, \ldots,k\}$. If the answer to all of them is \nop, then the answer to $\BrCPe(u,v,\varphi)$ is \nop. Otherwise, the answer to $\BrCPe(u,v,\varphi)$ is \yep.

It remains to consider the case where $\mathfrak{A}(u,v,\varphi)$ outputs $k$ (and hence $u \varphi^{k+1} \conj v \varphi$) and $\mathfrak{A}(v\varphi,v,\varphi)$ outputs $p$. Note that then $v\varphi^{p+1} \conj v \ker(\varphi)$ and hence $v \varphi^{p+2} \conj v\varphi$.
Moreover, 
\[u\varphi^{k+1} \sim v\varphi  \sim v\varphi^{p+2}=  (v\varphi)\varphi^{p+1}  \sim (u\varphi^{k+1})\varphi^{p+1}=u\varphi^{k+p+2},\]
which means that there are only finitely many conjugacy classes in the $\varphi$-orbit of $u$ and, again, it is enough to 
check whether $u \varphi^i \conj v$ for $i\in \{0,\ldots,k+p+1\}$. The answer to  $\BrCPe(u,v,\varphi)$ is \nop\ if the answer to all of them is \nop, and \yep\ otherwise.
This concludes the proof of the decidability of $\BrCPe(\Fn)$. 
\end{proof}

Finally, we  note that this work has continuation in \cite{carvalho_orbit_2024} and \cite{carvalho_free-abelian_2024}, where the same authors study Brinkmann Problems in some natural extensions of free groups.

\section*{Acknowledgements}
The first author is funded by national funds through the FCT - Funda\c c\~ao para a Ci\^encia e a
Tecnologia, I.P., under the scope of the projects UIDB/00297/2020 and UIDP/00297/2020 (Center for Mathematics and Applications).

The second author acknowledges support from the Spanish Agencia Estatal de Investigación
through grant PID2021-126851NB-I00 (AEI/ FEDER, UE), as well as from the Universitat Politècnica
de Catalunya in the form of a ``María Zambrano'' scholarship.



\renewcommand*{\bibfont}{\small}
\printbibliography

\Addresses

\end{document}

